\documentclass[12pt,a4paper]{article}
\usepackage[utf8]{inputenc}
\usepackage{amsfonts,amsmath,amssymb,amsthm}
\usepackage{hyperref}

\newtheorem{theorem}{Theorem}

\newtheorem{corollary}{Corollary}

\theoremstyle{remark}
\newtheorem{remark}{Remark}

\begin{document}

\title{Analytical and easily calculated expressions for continuous commutation functions under Gompertz-Makeham mortality}
\author{Andreas Nordvall Lager{\aa}s\footnote{Address: Department of Mathematics, Stockholm University, SE-106 91 Stockholm, Sweden. E-mail:\href{mailto:andreas@math,su,se}{\texttt{andreas@math.su.se}}}}
\date{}
\maketitle

\begin{abstract}
\noindent It is known, but perhaps not well-known, that when the mortality is assumed to be of Gompertz-Makeham-type, the expected remaining life-length and the commutation functions used for calculating the expected values of various types of life insurances can be expressed with an incomplete gamma function with a negative shape parameter. This is not of much use if ones software cannot calculate these values. The aim of this note is to show that one can express the commutation functions using only the exponential function, the (ordinary) gamma function and the gamma distribution function, which are all implemented in common statistical and spreadsheet software. This eliminates the need to evaluate the commutation functions and expected remaining life-length with numerical integration.\\

\end{abstract}

We will assume that individual lifes are distributed according to the Gomp\-ertz-Makeham distribution with distribution function $F(x;\alpha,\beta,\gamma)=1-l(x;\alpha,\beta,\gamma)$ and survival function $l(x;\alpha,\beta,\gamma)=\exp\{-\alpha x-\frac{\beta}{\gamma}(e^{\gamma x}-1)\}$ corresponding to the mortality rate $\mu(x;\alpha,\beta,\gamma)=-\frac{d}{dx}l(x;\alpha,\beta,\gamma)=\alpha + \beta e^{\gamma x}$ at age $x$. We assume a fixed continuously compounded interest rate $\delta$. Let $\bar{a}_x=\bar{a}_x(\alpha,\beta,\gamma,\delta)$ be the present value of a life-long annuity, continuously paid at rate 1, to a person at age $x$, and $e_x(\alpha,\beta,\gamma)$ the expected remaining life-length for a person at age $x$. The main result is the following.
\begin{theorem}
\begin{align}
\bar{a}_x(\alpha,\beta,\gamma,\delta) &= e_0(\alpha+\delta,\beta e^{\gamma x},\gamma),\text{where}\label{ax}\\
e_0(\alpha,\beta,\gamma)&=\tfrac{1}{\alpha}\Big(1-\big(\tfrac{\beta}{\gamma}\big)^{\alpha/\gamma}e^{\beta/\gamma}\Gamma(1-\tfrac{\alpha}{\gamma})\big[1-G(\tfrac{\beta}{\gamma};1-\tfrac{\alpha}{\gamma},1)\big]\Big),\text{i.e.}\label{e0}\\
\bar{a}_x(\alpha,\beta,\gamma,\delta)&=\tfrac{1}{\alpha+\delta}\Big(1-\big(\tfrac{\beta e^{\gamma x}}{\gamma}\big)^{(\alpha+\delta)/\gamma}e^{\beta e^{\gamma x}/\gamma}\times\notag\\
&\qquad\qquad\qquad\times\Gamma(1-\tfrac{\alpha+\delta}{\gamma})\big[1-G(\tfrac{\beta e^{\gamma x}}{\gamma};1-\tfrac{\alpha+\delta}{\gamma},1)\big]\Big),\label{ax_hel}
\end{align}
where $\Gamma(\eta)$ is the gamma function: $\Gamma(\eta)=\int_0^{\infty}y^{\eta-1}e^{-y}dy$, and $G(z;\eta,1)$ is the distribution function of a gamma distributed random variable with shape parameter $\eta$ and scale parameter 1: $G(z;\eta,1)=\frac{1}{\Gamma(\eta)}\int_0^z y^{\eta-1}e^{-y}dy$.
\end{theorem}
\begin{remark}
Even though equation \eqref{ax_hel} may look unwieldy, it only consists of functions that are readily available in common statistical and spreadsheet software, e.g.\ in R and Microsoft Excel. Thus, \emph{no numerical integration is necessary} to evaluate these expressions. Note, however, the findings of Yalta (2008).
\end{remark}
\begin{remark}\label{ageing}
The special case with zero mortality, i.e.\ $\alpha=\beta=0$ (and $\gamma$ arbitrary), yields $\bar{a}_x(0,0,\gamma,\delta) = e_0(\delta,0,\gamma)=1/\delta$ as expected for a continuously compounded perpetuity. With $\alpha>0$ and $\beta=0$, the life-length is exponentially distributed and there is no ageing. Then $\bar{a}_x(\alpha,0,\gamma,\delta) = e_0(\alpha+\delta,0,\gamma)=1/(\alpha+\delta)$. The general case has $\bar{a}_x(\alpha,\beta,\gamma,\delta)=\tfrac{1}{\alpha+\delta}(1-[\dots])$ and one can thus interpret the ellipsis part $[\dots]$ as the \emph{effect of ageing} on the value of the annuity.
\end{remark}

With $\delta=0$, the present value of the annuity equals the nominal value, and since the annuity is paid continuously at rate 1, this nominal value is equal to the remaining life-length. Thus, we arrive at the following result, which will also be reached through another route in the proof of Theorem 1.
\begin{corollary}\label{ex}
$e_x(\alpha,\beta,\gamma)=\bar{a}_x(\alpha,\beta,\gamma,0)=e_0(\alpha,\beta e^{\gamma x},\gamma)$.
\end{corollary}

In order to calculate the present values of different life insurance contracts one usually defines the commutation functions
\begin{align*}
D(x) &= l(x)e^{-\delta x},\\
N(x) &= \int_x^{\infty}D(y)dy,\\
M(x) &= \int_x^{\infty}\mu(y)D(y)dy.
\end{align*}
The corresponding commutation functions with twice the interest rate are also useful when one wants to calculate not only the expected value of different insurances, but also their variance, see Andersson (2005). In general $M(x)=D(x)-\delta N(x)$, and with Gompertz-Makeham mortality $D(x)=l(x;\alpha+\delta,\beta,\gamma)$, which is easy to evaluate numerically. The only ``difficult'' function is $N(x)$. However, a standard result is that $\bar{a}_x=N(x)/D(x)$, so by Theorem 1 we get a tractable expression also for $N(x)$:
\begin{corollary}\label{Nx}
$N(x) = D(x)e_0(\alpha+\delta,\beta e^{\gamma x},\gamma)$, where $D(x)=l(x;\alpha+\delta,\beta,\gamma)$.
\end{corollary}

It remains to prove Theorem 1.

\begin{proof}[Proof of Theorem 1] The incomplete gamma function is defined as 
$
\Gamma(\eta,z)=\int_z^{\infty}y^{\eta-1}e^{-y}dy.
$ 
In analogy with the gamma distribution we call $\eta$ the shape parameter. The incomplete gamma function is closely related to the gamma distribution function:
\begin{align*}
G(z;\eta,1)&=\frac{1}{\Gamma(\eta)}\int_0^z y^{\eta-1}e^{-y}dy=1-\frac{1}{\Gamma(\eta)}\int_z^{\infty}y^{\eta-1}e^{-y}dy = 1 - \frac{\Gamma(\eta,z)}{\Gamma(\eta)},
\end{align*}
or, equivalently,
\begin{equation}\label{gam_gam}
\Gamma(\eta,z)=\Gamma(\eta)(1-G(z;\eta,1)),
\end{equation}
for $\eta>0$. By partial integration,
\begin{align*}
\Gamma(\eta,z)&=\int_z^{\infty}y^{\eta-1}e^{-y}dy=\frac{1}{\eta}\left[y^{\eta}e^{-y}\right]_z^{\infty}+\frac{1}{\eta}\int_z^{\infty}y^{\eta}e^{-y}dy\\ 
&=\frac{1}{\eta}\left(\Gamma(\eta+1,z)-z^{\eta}e^{-z}\right). 
\end{align*}
Combining this with \eqref{gam_gam} yields
\begin{equation}\label{gamma_dist}
\Gamma(\eta,z)=-\frac{1}{\eta}\big(z^{\eta}e^{-z}-\Gamma(\eta+1)\big[1-G(z;\eta+1,1)\big]\big),
\end{equation}
which we will use shortly.

Following Andersson (2005), we express $e_0$ with an incomplete gamma function:
\begin{align}
e_0(\alpha,\beta,\gamma)&=\int_0^{\infty}l(t;\alpha,\beta,\gamma)dt = \int_0^{\infty}e^{-\alpha t-\frac{\beta}{\gamma}(e^{\gamma t}-1)}dt \label{e0_int}\\
\left\{y=\frac{\beta}{\gamma}e^{\gamma t}\right\}&= \left(\frac{\beta}{\gamma}\right)^{\alpha/\gamma}\frac{e^{\beta/\gamma}}{\gamma}\int_{\beta/\gamma}^{\infty}y^{-\frac{\alpha}{\gamma}-1}e^{-y}dy\notag\\
&=\left(\frac{\beta}{\gamma}\right)^{\alpha/\gamma}\frac{e^{\beta/\gamma}}{\gamma}\Gamma\left(-\frac{\alpha}{\gamma},\frac{\beta}{\gamma}\right)\label{e0_gamma}.
\end{align}

We continue with the expected remaining life-length at age $x$:
$$
e_x(\alpha,\beta,\gamma)=\int_0^{\infty}\frac{l(x+t;\alpha,\beta,\gamma)}{l(x;\alpha,\beta,\gamma)}dt = \int_0^{\infty}e^{-\alpha t-\frac{\beta}{\gamma}e^{\gamma x}(e^{\gamma t}-1)}dt,
$$
and by comparing with \eqref{e0_int} we arrive at the conclusion of Corollary \ref{ex}: $e_x(\alpha,\beta,\gamma)=e_0(\alpha,\beta e^{\gamma x},\gamma)$.

By standard arguments,
$$
\bar{a}_x(\alpha,\beta,\gamma,\delta) = \int_0^{\infty}e^{-\delta t}\frac{l(x+t;\alpha,\beta,\gamma)}{l(x;\alpha,\beta,\gamma)}dt= \int_0^{\infty}e^{-(\alpha+\delta)t-\frac{\beta}{\gamma}e^{\gamma x}(e^{\gamma t}-1)}dt,
$$
and by once again comparing with \eqref{e0_int} equation \eqref{ax} in the theorem is proved.

Equation \eqref{e0} is obtained by using \eqref{gamma_dist} on \eqref{e0_gamma} to rewrite the incomplete gamma function in the latter equation. Equation \eqref{ax_hel} follows immediately from \eqref{ax} and \eqref{e0}. This concludes the proof of Theorem 1.
\end{proof} 

\begin{remark}
Note that for typical values for the parameters, e.g.\ $\alpha=0.001, \beta=0.000012, \gamma=0.101314$ and $\delta=0.026559$ from Andersson (2005), we have in equation \eqref{ax_hel} a shape parameter $1-(\alpha+\delta)/\gamma\doteq 0.727984$ \emph{which is positive}, just as shape parameters for gamma distributions should be. If one has to ones disposal software that can calculate incomplete gamma functions with negative shape parameters, then one can use equation \eqref{e0_gamma} straight away, see e.g.\ Mingari Scarpello et al.\ (2006). Since the effect of ageing on the value of an annuity, cf.\ Remark \ref{ageing}, is not easily seen from \eqref{e0_gamma}, the formulation in the theorem is nevertheless of some independent theoretical interest, regardless of whether it is easy to evaluate functions with negative shape parameters or not.
\end{remark}

\section*{References}

\textsc{Andersson, G.} (2005). \textit{Livf{\"o}rs{\"a}kringsmatematik}. Svenska F{\"o}r\-s{\"a}k\-rings\-f{\"o}ren\-ingen, Stockholm. ISBN: 91-974960-1-4

\textsc{Mingari Scarpello, G.; Ritelli, D. and Spelta, D.} (2006). Actuarial values calculated using the incomplete gamma function. \href{http://rivista-statistica.cib.unibo.it/issue/view/50}{\emph{Statistica (Bologna)}} \textbf{66}, no.\ 1, 77--84.

\textsc{Yalta, A.T.} (2008). The accuracy of statistical distributions in Microsoft\raisebox{1ex}{\scriptsize\textregistered} Excel. \textit{Computational Statistics and Data Analysis} \textbf{52}, 4579--4586.\\
\href{http://dx.doi.org/10.1016/j.csda.2008.03.005}{\texttt{doi:10.1016/j.csda.2008.03.005}}

\end{document}